\documentclass{amsart}
\usepackage{amssymb,amsthm,amsmath}
\usepackage[numbers,sort&compress]{natbib}

\thanks{V. R\u adulescu and D. Repov\v{s} were supported in part by  the  Slovenian  Research  Agency
grants P1-0292, N1-0064, J1-8131, J1-7025, and J1-6721. V. R\u adulescu acknowledges the support through a grant of the Romanian National
Authority for Scientific Research and Innovation, CNCS-UEFISCDI, project number PN-III-P4-ID-PCE-2016-0130. Q. Zhang has been partially supported by the key projects of Science and
Technology Research of the Henan Education Department (14A110011) and the
National Natural Science Foundation of China (11326161 and 10971087).}

\newcommand{\di}{\displaystyle}

\newcommand{\RR}{\mathbb{R}}

\theoremstyle{plain}

\theoremstyle{definition}

\begin{document}
\title[Nonhomogeneous problems without Ambrosetti-Rabinowitz condition]{
Nonhomogeneous Dirichlet problems without the Ambrosetti-Rabinowitz condition
}
\author[G. Li, V. D. R\u{a}dulescu, D. D. Repov\v{s} and
Q. Zhang]{Gang Li, Vicen\c{t}iu D. R\u{a}dulescu, Du\v{s}an D. Repov\v{s} and
Qihu Zhang}
\address[G. Li]{College of Mathematics and Information Science, Zhengzhou
University of Light Industry, Zhengzhou, Henan 450002, China}
\email{leagongpaper@yeah.net}
\address[V.D. R\u{a}dulescu]{Department of Mathematics, Faculty of Sciences,
King Abdulaziz University, P.O. Box 80203, Jeddah 21589, Saudi Arabia \&
Department of Mathematics, University of Craiova, 200585 Craiova, Romania}
\email{vicentiu.radulescu@math.cnrs.fr}
\address[D.D. Repov\v{s}]{Faculty of Mathematics
and Physics, University of Ljubljana, Jadranska 19, SI-1000
Ljubljana, Slovenia}
\email{dusan.repovs@guest.arnes.si}
\address[Q. Zhang]{College of Mathematics and Information Science, Zhengzhou
University of Light Industry, Zhengzhou, Henan 450002, China}
\email{zhangqihu@yahoo.com}

\begin{abstract}
We consider the existence of solutions of the following $p(x)$-Laplacian
Dirichlet problem without the Ambrosetti-Rabinowitz condition:
\begin{equation*}
\left\{
\begin{array}{l}
-\mbox{div}(|\nabla u|^{p(x)-2}\nabla u)=f(x,u),\text{ in }\Omega , \\
u=0,\text{ on }\partial \Omega .%
\end{array}
\right.
\end{equation*}
We give a new growth condition and we point out its importance for checking
the Cerami compactness condition. We prove the existence of solutions of the
above problem via the critical point theory, and also provide some
multiplicity properties. Our results extend previous results of Q.~Zhang and
C.~Zhao (Existence of strong solutions of a $p(x)$-Laplacian Dirichlet
problem without the Ambrosetti-Rabinowitz condition, \textit{Computers and
Mathematics with Applications}, 2015) and we establish the existence of
solutions under weaker hypotheses on the nonlinear term.
\end{abstract}

\maketitle

%%% ----------------------------------------------------------------------

\baselineskip 14pt

\section{ Introduction}

In recent years, the study of differential equations and variational problems
with variable exponent growth conditions has been a topic of great interest.
This type of problems has very strong background, for instance in image
processing, nonlinear electro-rheological fluids and elastic mechanics. Some
of these phenomena are related to the \textit{Winslow effect}, which
describes the behavior of certain fluids that become solids or quasi-solids
when subjected to an electric field. The result was named after the American
engineer Willis M. Winslow.

There are many papers dealing with problems with variable exponents, see
\cite{e1}-\cite{j8}, \cite{e2}-\cite{e20}, \cite{e21}, \cite{a1}-\cite%
{radnla}, \cite{e3}, \cite{e26}, \cite{e27}-\cite{55}, \cite{e4}-\cite{e34}.
On the existence of solutions of these kinds of problems, we refer to \cite%
{j8, e11, e12, e16, e17, a1, a3, e32}. We also refer to the recent monograph
\cite{radrep} dealing with variational methods in the framework of nonlinear
problems with variable exponent.

In this paper, we consider the existence of solutions of the following class
of Dirichlet problems:
\begin{equation*}
\text{(P) }\ \left\{
\begin{array}{l}
-\Delta _{p(x)}u:=-\mbox{div}\,(\left\vert \nabla u\right\vert
^{p(x)-2}\nabla u)=f(x,u),\text{ }\left. {}\right. \text{ in }\Omega , \\
u=0,\text{ }\left. {}\right. \text{on }\partial \Omega ,%
\end{array}
\right.
\end{equation*}
where $\Omega \subset \mathbb{R}^{N}$ is a bounded domain with $C^{1,\alpha
} $ smooth boundary, and $p(\cdot )>1$ is of class $C^{1}(\overline{ \Omega }%
)$.

Since the elliptic operator with variable exponent is not homogeneous, new
methods and techniques are needed to study these types of problems. We point
out that commonly known methods and techniques for studying constant
exponent equations fail in the setting of problems involving variable
exponents. For instance, the eigenvalues of the $p(x)$-Laplacian Dirichlet
problem were studied in \cite{e13}. In this case, if $\Omega \subset \mathbb{%
R}^{N}$ is a smooth bounded domain, then the Rayleigh quotient
\begin{equation}
\lambda _{p(\cdot )}=\displaystyle\underset{u\in W_{0}^{1,p(\cdot )}(\Omega
)\backslash \{0\}}{\inf }\frac{\di\int_{\Omega }\frac{1}{p(x)}\left\vert \nabla
u\right\vert ^{p(x)}dx}{\di\int_{\Omega }\frac{1}{p(x)}\left\vert u\right\vert
^{p(x)}dx}  \label{d2}
\end{equation}%
is in general zero, and $\lambda _{p(\cdot )}>0$ holds only under some
special conditions.

In \cite{j7}, the author generalized the Picone identities for half-linear
elliptic operators with $p(x)$-Laplacian. In the same paper some
applications to Sturmian comparison theory are also presented, but the
formula is different from the constant exponent case. In a related setting,
we point out that the formula
\begin{equation*}
\int_\Omega |u(x)|^pdx=p\int_0^\infty t^{p-1}\,|\{x\in\Omega ;\
|u(x)|>t\}|\,dt
\end{equation*}
has no variable exponent analogue.

In \cite{j6} and \cite{55} the authors deal with the local boundedness and
the Harnack inequality for the $p(x)$-Laplace equation. But it was shown in
\cite{j6} that even in the case of a very nice exponent, for example,
\begin{equation*}
p(x):=\displaystyle\left\{
\begin{array}{lll}
& \displaystyle 3, & \quad\text{ for }0<x\leq\displaystyle \frac{1}{2} \\
& \displaystyle 3-2\left(x-\frac{1}{2}\right), & \quad\text{ for }%
\displaystyle\frac{1}{2}<x<1%
\end{array}%
\right.
\end{equation*}%
the constant in the Harnack inequality depends on the minimizer, that is,
the inequality $\sup u\leq c\inf u$ does not hold for any absolute constant $%
c$.

%Variable exponent Lebesgue spaces do not have the \textit{mean continuity
%property}: if $p$ is continuous and nonconstant in an open ball $B$, then
%there exists a function $u\in L^{p(\cdot )}(B)$ such that $u(x+h)\not\in
%L^{p(\cdot )}(B)$ for all $h\in \mathbb{R}^{N}$ with an arbitrary small norm.

The standard norm in variable exponent Sobolev spaces is the so-called
Luxemburg norm $\left\vert u\right\vert _{p(\cdot )}$ (see section 2) and
the integral $\int_{\Omega }\left\vert u(x)\right\vert ^{p(x)}dx$ does not
satisfy the constant power relation.

In many instances, it is difficult to judge whether or not results about $%
p$-Laplacian can be generalized to $p(x)$-Laplacian, and even if this can be
done, it is still difficult to figure out the form in which the results
should be.

Our main goal is to obtain a couple of existence results for the problem (P)
without the Ambrosetti-Rabinowitz condition via critical point theory. For
this purpose, we use a new method for checking the Cerami compactness
condition under a new growth condition. Our results can be regarded as
extensions of the corresponding results for the $p$-Laplacian problems, but
the growth condition and the methods for checking the Cerami compactness
condition are different with respect to quasilinear equations with constant
exponent.

Next, we give a review of some results related to our work. Since the
Ambrosetti-Rabinowitz type condition is quite restrictive and excludes many
cases of nonlinearity, there are many papers dealing with the problem
without the Ambrosetti-Rabinowitz type growth condition. For the constant
exponent case $p(\cdot )\equiv p$, we refer to \cite{jj5, jj3, jj2, jj4}.

In \cite{jj5}, the authors considered the problem (P) for $p(\cdot )\equiv p
$, and proved the existence of weak solutions under the following
assumptions: $\underset{\left\vert t\right\vert \rightarrow +\infty }{\lim }$
$\frac{F(x,t)}{\left\vert t\right\vert ^{p}}=+\infty $, where $%
F(x,t)=\int_{0}^{t}f(x,s)ds$; and there exists a constant $C_{\ast }>0$ such
that $H(x,t)\leq H(x,s)+C_{\ast }$ for each $x\in \Omega $, $0<t<s$ or $%
s<t<0 $, where $H(x,t)=$ $tf(x,t)-pF(x,t)$.

In \cite{jj3}, the author studied the problem (P) for $p(\cdot )\equiv p$.
Under the assumption that $\frac{f(x,s)}{\left\vert s\right\vert ^{p-2}s}$
is increasing when $s\geq s_{0}$ and decreasing when $s\leq -s_{0}$, $%
\forall x\in \Omega $, the existence of weak solutions was obtained.

In \cite{jj2}, the authors studied the problem (P) for $p(\cdot )\equiv 2$,
which becomes a Laplacian problem. The main result in \cite{jj2} establishes
the existence of weak solutions by assuming that $\frac{f(x,s)}{s}$ is
increasing when $s\geq s_{0}$ and decreasing when $s\leq -s_{0}$, for all $%
x\in \Omega $.

In \cite{jj4}, the author also studied the problem (P) for $p(\cdot )\equiv
2 $ and proved the existence of weak solutions under the assumption
\begin{equation*}
sf(x,u)\geq C_{0}\left\vert s\right\vert ^{\mu }\text{, where }\mu >2\text{
and }C_{0}>0.
\end{equation*}

If $p(\cdot )$ is a general function, results on variable exponent problem
without the Ambrosetti-Rabinowitz type growth condition are rare due to the
complexity of $p(x)$-Laplacian (see \cite{jj1, aa3, aa4, aa8, 31a}). However
their assumptions imply $G_{p^{+}}(x,t)=f(x,t)t-p^{+}F(x,t)\geq 0$ and $%
F(x,t)>0$ as $t\rightarrow +\infty $, so we can see that $F(x,t)\geq
Ct^{p^{+}}$ as $t\rightarrow +\infty $. This is too strong and unnatural for
the $p(x)$-Laplacian problems.

In \cite{e32}, the author considered the problem (P) under the following
growth condition: \newline
there exist constants $M,C_{1},C_{2}>0,a>p$ on $\overline{\Omega }$ such
that
\begin{equation}
C_{1}\left\vert t\right\vert ^{p(x)}[\ln (e+\left\vert t\right\vert
)]^{a(x)-1}\leq C_{2}\frac{tf(x,t)}{\ln (e+\left\vert t\right\vert )}\leq
tf(x,t)-p(x)F(x,t),\forall \left\vert t\right\vert \geq M,\forall x\in
\Omega .  \label{ja1}
\end{equation}

A typical example is $f(x,t)=\left\vert t\right\vert ^{p(x)-2}t[\ln
(1+\left\vert t\right\vert )]^{a(x)}$. This function satisfies the above
condition (\ref{ja1}), but does not satisfy the Ambrosetti-Rabinowitz
condition.

Our paper was motivated by \cite{e32}. We further weaken the condition (\ref%
{ja1}). To begin we point out that the assumption $a>p$ on $\overline{\Omega
}$ is unnecessary in the present paper.

Before stating our main results, we make the following assumptions:

(f$_{0}$): $f:\Omega \times \mathbb{R}\rightarrow \mathbb{R}$ satisfies the
Carath\'eodory condition and
\begin{equation*}
\left\vert f(x,t)\right\vert \leq C(1+\left\vert t\right\vert ^{\alpha
(x)-1}),\forall (x,t)\in \Omega \times \mathbb{R},
\end{equation*}%
where $\alpha \in C(\overline{\Omega })$ and $p(x)<\alpha (x)<p^{\ast }(x)$
on $\overline{\Omega }$.

(f$_{1}$): there exist constants $M,C>0$,$\ $such that%
\begin{equation}
C\frac{tf(x,t)}{K(t)}\leq tf(x,t)-p(x)F(x,t),\forall \left\vert t\right\vert
\geq M,\forall x\in \overline{\Omega },  \label{ja2}
\end{equation}%
and
\begin{equation}
\frac{tf(x,t)}{\left\vert t\right\vert ^{p(x)}[K(t)]^{p(x)}}\rightarrow
+\infty \text{ uniformly as }\left\vert t\right\vert \rightarrow +\infty
\text{ for }x\in \overline{\Omega },  \label{ja3}
\end{equation}%
where $K$ satisfies the following hypotheses:

(K): $1\leq K(\cdot )\in C^{1}([0,+\infty ),[1,+\infty ))$ is increasing and
$[\ln (e+t)]^{2}\geq K(t)\rightarrow +\infty $ as $\left\vert t\right\vert
\rightarrow +\infty $, which satisfies $tK^{\prime }(t)/K(t)\leq \sigma
_{0}\in (0,1)$, where $\sigma _{0}$ is a constant.

(f$_{2}$): $f(x,t)=o(\left\vert t\right\vert ^{p(x)-1})$ uniformly for $x\in
\Omega $ as $t\rightarrow 0$.

(f$_{3}$): $f(x,-t)=-f(x,t),$ $\forall x\in \overline{\Omega },$ $\forall
t\in \mathbb{R}$.

(f$_{4}$): $F$ satisfies
\begin{equation*}
\frac{F(x,t)}{\left\vert t\right\vert ^{p(x)}[\ln (e+\left\vert t\right\vert
)]^{p(x)}}\rightarrow +\infty \text{ uniformly as }\left\vert t\right\vert
\rightarrow +\infty \text{ for }x\in \overline{\Omega }.
\end{equation*}

(p$_{1}$): there is a vector $l\in \mathbb{R}^{N}\backslash \{0\}$ such that
for any $x\in \Omega $, $\rho (t)=p(x+tl)$ is monotone for $t\in
I_{x}(l)=\{t\mid x+tl\in \Omega \}$.

(p$_{2}$): $p$ has a local maximum point, that is, there exist $x_{0}\in
\Omega $ and $\delta >0$ such that $\overline{B(x_{0},3\delta )}\subset
\Omega $ and
\begin{equation*}
\underset{\left\vert x-x_{0}\right\vert \leq \delta }{\min }p(x)>\underset{%
2\delta \leq \left\vert x-x_{0}\right\vert \leq 3\delta }{\max }p(x).
\end{equation*}

(p$_{3}$): $p$ has a sequence of local maximum points, that is, there exist
a sequence of points $x_{n}\in \Omega $ and $\delta _{n}>0$ such that $%
\overline{B(x_{0},3\delta _{n})}$ are mutually disjoint and
\begin{equation*}
\underset{\left\vert x-x_{n}\right\vert \leq \delta _{n}}{\min }p(x)>%
\underset{2\delta _{n}\leq \left\vert x-x_{n}\right\vert \leq 3\delta _{n}}{%
\max }p(x).
\end{equation*}

We state our main results in what follows.

\medskip \textbf{Theorem 1.1 } Assume that hypotheses (f$_{0}$)-(f$_{2}$), (p%
$_{1}$), and (f$_{4}$) or (p$_{2}$) are fulfilled. Then problem (P) has a
nontrivial solution.\newline

\textbf{Theorem 1.2 } Assume that hypotheses (f$_{0}$), (f$_{1}$), (f$_{3}$%
), and (f$_{4}$) or (p$_{3}$) are fulfilled. Then problem (P) has infinitely
many pairs of solutions.\newline

\textbf{Remark}. (i). The following functions satisfy the hypothesis ($K$):
\begin{eqnarray*}
K_{1}(t) &=&\ln (e+\left\vert t\right\vert ) \\
K_{2}(t) &=&\ln (e+\ln (e+\left\vert t\right\vert )) \\
K_{3}(t) &=&[\ln (e+\ln (e+\left\vert t\right\vert ))]\ln (e+\left\vert
t\right\vert ).
\end{eqnarray*}%
Let $K=K_{1}$, and $f(x,t)=\left\vert t\right\vert ^{p(x)-2}t[\ln
(1+\left\vert t\right\vert )]^{p(x)}\rho (\left\vert t\right\vert )$, where $%
1\leq \rho (\left\vert t\right\vert )\leq \lbrack \ln (e+\left\vert
t\right\vert )]^{2}$, $\rho ^{\prime }\geq 0$ and $\rho (\left\vert
t\right\vert )\rightarrow +\infty $ as $\left\vert t\right\vert \rightarrow
+\infty $, for example $\rho (\left\vert t\right\vert )=\ln (e+\ln
(e+\left\vert t\right\vert ))$. Then $f$ satisfies the condition (f$_{0}$)-(f%
$_{4}$), but it does not satisfy the Ambrosetti-Rabinowitz condition, and
does not satisfy (\ref{ja1}) (ii). We do not need any monotonicity
assumption on $f(x,\cdot )$.

This paper is organized as follows. In Section 2, we do some preparatory
work including some basic properties of the variable exponent Sobolev spaces,
which can be regarded as a special class of generalized Orlicz-Sobolev
spaces. In Section 3, we give proofs of the results stated above.

\section{ Preliminary results}

Throughout this paper, we use letters $c,\ c_{i},\ C,\ C_{i}$, $i=1,2,...$
to denote generic positive constants which may vary from line to line, and
we will specify them whenever necessary.

One of the reasons for the huge development of the theory of classical
Lebesgue and Sobolev spaces $L^{p}$ and $W^{1,p}$ (where $1\leq p\leq \infty
$) is its usefulness for the description of many phenomena arising in
applied sciences. For instance, many materials can be modeled with
sufficient accuracy by using the function spaces $L^{p}$ and $W^{1,p}$,
where $p$ is a fixed constant. For some materials with nonhomogeneities, for
instance electro-rheological fluids (sometimes referred to as
\textquotedblleft smart fluids"), this approach is not adequate, but rather
the exponent $p$ should be allowed to vary. This leads us to the study of
variable exponent Lebesgue and Sobolev spaces, $L^{p(\cdot )}$ and $%
W^{1,p(\cdot )}$, where $p$ is a real--valued function.

In order to discuss problem (P), we need some results about the space $%
W_{0}^{1,p(\cdot )}(\Omega )$, which we call \textit{variable exponent
Sobolev space}. We first state some basic properties of $W_{0}^{1,p(\cdot
)}(\Omega )$ (for details, see \cite{j5, e9, e12, e20, radrep, e26}). Denote
\begin{eqnarray*}
C_{+}(\overline{\Omega }) &=&\left\{ h\left\vert h\in C(\overline{\Omega })%
\text{, }h(x)>1\text{ for }x\in \overline{\Omega }\right. \right\} , \\
h^{+} &=&\max_{\overline{\Omega }}h(x)\text{, }h^{-}=\min_{\overline{\Omega }%
}h(x)\text{, for any }h\in C(\overline{\Omega }), \\
L^{p(\cdot )}(\Omega ) &=&\left\{ u\mid u\text{ is a measurable real-valued
function, }\int_{\Omega }\left\vert u(x)\right\vert ^{p(x)}dx<\infty
\right\} .
\end{eqnarray*}

We introduce the norm on $L^{p(\cdot )}(\Omega )$ by

\begin{equation*}
\left\vert u\right\vert _{p(\cdot )}=\inf \left\{ \lambda >0\left\vert
\int_{\Omega }\left\vert \frac{u(x)}{\lambda }\right\vert ^{p(x)}dx\leq
1\right. \right\}\,.
\end{equation*}

Then ($L^{p(\cdot )}(\Omega )$, $\left\vert \cdot \right\vert _{p(\cdot )}$)
becomes a Banach space and it is called the \textit{variable exponent
Lebesgue space.} \newline

\textbf{Proposition 2.1 }(see \cite{j5, radrep}). i) The space $(L^{p(\cdot
)}(\Omega ),\left\vert \cdot \right\vert _{p(\cdot )})$ is a separable,
uniform convex Banach space, and its conjugate space is $L^{q(\cdot
)}(\Omega ),$ where $\frac{1}{q(\cdot )}+\frac{1}{p(\cdot )}\equiv 1$. For
any $u\in L^{p(\cdot )}(\Omega )$ and $v\in L^{q(\cdot )}(\Omega )$, we have
\begin{equation*}
\left\vert \int_{\Omega }uvdx\right\vert \leq \left( \frac{1}{p^{-}}+\frac{1%
}{q^{-}}\right) \left\vert u\right\vert _{p(\cdot )}\left\vert v\right\vert
_{q(\cdot )}.
\end{equation*}

ii) If $p_{1},$ $p_{2}\in C_{+}(\overline{\Omega })$ , $p_{1}(x)\leq
p_{2}(x) $ for any $x\in \overline{\Omega },$ then $L^{p_{2}(\cdot )}(\Omega
)\subset L^{p_{1}(\cdot )}(\Omega ),$ and this imbedding is continuous.%
\newline

\textbf{Proposition 2.2 } (see \cite{e12,radrep}). If $f:$ $\Omega \times
%TCIMACRO{\U{211d} }%
%BeginExpansion
\mathbb{R}
%EndExpansion
\rightarrow
%TCIMACRO{\U{211d} }%
%BeginExpansion
\mathbb{R}
%EndExpansion
$ is a Carath\'eodory function and satisfies%
\begin{equation*}
\left\vert f(x,s)\right\vert \leq a(x)+b\left\vert s\right\vert
^{p_{1}(x)/p_{2}(x)}\left. {}\right. \text{for any }x\in \Omega ,s\in
%TCIMACRO{\U{211d} }%
%BeginExpansion
\mathbb{R}
%EndExpansion
,
\end{equation*}%
where $p_{1}$, $p_{2}\in C_{+}(\overline{\Omega })$ , $a\in L^{p_{2}(\cdot
)}(\Omega )$, $a(x)\geq 0$, $b\geq 0$, then the Nemytsky operator from $%
L^{p_{1}(\cdot )}(\Omega )$ to $L^{p_{2}(\cdot )}(\Omega )$ defined by $%
(N_{f}u)(x)=f(x,u(x))$, is a continuous and bounded operator. \newline

\textbf{Proposition 2.3 }(see \cite{e12,radrep}). If we denote
\begin{equation*}
\rho (u)=\int_{\Omega }\left\vert u\right\vert ^{p(x)}dx\text{, }\forall
u\in L^{p(\cdot )}(\Omega ),
\end{equation*}%
then there exists $\xi \in \overline{\Omega }$ such that $\left\vert
u\right\vert _{p(\cdot )}^{p(\xi )}=\int_{\Omega }\left\vert u\right\vert
^{p(x)}dx$ and

i) $\left\vert u\right\vert _{p(\cdot )}<1(=1;>1)\Longleftrightarrow \rho
(u)<1(=1;>1);$

ii) $\left\vert u\right\vert _{p(\cdot )}>1\Longrightarrow \left\vert
u\right\vert _{p(\cdot )}^{p^{-}}\leq \rho (u)\leq \left\vert u\right\vert
_{p(\cdot )}^{p^{+}};$ $\left\vert u\right\vert _{p(\cdot
)}<1\Longrightarrow \left\vert u\right\vert _{p(\cdot )}^{p^{-}}\geq \rho
(u)\geq \left\vert u\right\vert _{p(\cdot )}^{p^{+}};$

iii) $\left\vert u\right\vert _{p(\cdot )}\rightarrow 0\Longleftrightarrow
\rho (u)\rightarrow 0;$ $\left\vert u\right\vert _{p(\cdot )}\rightarrow
\infty \Longleftrightarrow \rho (u)\rightarrow \infty .$\newline

\textbf{Proposition 2.4 }(see \cite{e12,radrep}). If $u$, $u_{n}\in
L^{p(\cdot )}(\Omega )$, $n=1,2,\cdots ,$ then the following statements are
equivalent:

1) $\underset{k\rightarrow \infty }{\lim }$ $\left\vert u_{k}-u\right\vert
_{p(\cdot )}=0;$

2) $\underset{k\rightarrow \infty }{\lim }$ $\rho \left( u_{k}-u\right) =0;$

3) $u_{k}$ $\rightarrow $ $u$ in measure in $\Omega $ and $\underset{%
k\rightarrow \infty }{\lim }$ $\rho \left( u_{k}\right) =\rho (u).$\newline

The space $W^{1,p(\cdot )}(\Omega )$ is defined by%
\begin{equation*}
W^{1,p(\cdot )}(\Omega )=\left\{ u\in L^{p(\cdot )}\left( \Omega \right)
\mid \nabla u\in (L^{p(\cdot )}\left( \Omega \right) )^{N}\right\} ,
\end{equation*}%
and it can be equipped with the norm
\begin{equation*}
\left\Vert u\right\Vert =\left\vert u\right\vert _{p(\cdot )}+\left\vert
\nabla u\right\vert _{p(\cdot )},\forall u\in W^{1,p(\cdot )}\left( \Omega
\right) .
\end{equation*}

We denote by $W_{0}^{1,p(\cdot )}(\Omega )$ the closure of $C_{0}^{\infty
}\left( \Omega \right) $ in $W^{1,p(\cdot )}(\Omega )$ and set
\begin{equation*}
p^{\ast }(x)=\left\{
\begin{array}{l}
Np(x)/(N-p(x))\text{ , }p(x)<N, \\
\infty \text{ , }p(x)\geq N.%
\end{array}%
\right.
\end{equation*}

Then we have the following properties.\newline

\textbf{Proposition 2.5 }(see \cite{j5,e12,radrep}). i) $W^{1,p(\cdot
)}(\Omega )$ and $W_{0}^{1,p(\cdot )}(\Omega )$ are separable reflexive
Banach spaces;

ii) if $q\in C_{+}\left( \overline{\Omega }\right) $ and $q(x)<p^{\ast }(x)$
for any $x\in \overline{\Omega },$ then the imbedding from $W^{1,p(\cdot
)}(\Omega )$ to $L^{q(\cdot )}\left( \Omega \right) $ is compact;

iii) there is a constant $C>0$ such that
\begin{equation*}
\left\vert u\right\vert _{p(\cdot )}\leq C\left\vert \nabla u\right\vert
_{p(\cdot )},\forall u\in W_{0}^{1,p(\cdot )}(\Omega ).
\end{equation*}

It follows from iii) of Proposition 2.5 that $\left\vert \nabla u\right\vert
_{p(\cdot )}$ and $\left\Vert u\right\Vert $ are equivalent norms on $%
W_{0}^{1,p(\cdot )}(\Omega )$. From now on, we will use $\left\vert \nabla
u\right\vert _{p(\cdot )}$ instead of $\left\Vert u\right\Vert $ as the norm
on $W_{0}^{1,p(\cdot )}(\Omega )$. 

The  Lebesgue and Sobolev
spaces with variable exponents coincide with the usual Lebesgue and Sobolev spaces provided that $p$ is constant. These function spaces $L^{p(x)}$ and $W^{1,p(x)}$ have some non-usual properties, see \cite[p. 8-9]{radrep}. Some of these properties are the following:

(i)
Assuming that $1<p^-\leq p^+<\infty$ and $p:\overline\Omega\rightarrow [1,\infty)$ is a smooth function, then the following co-area formula
$$\int_\Omega |u(x)|^pdx=p\int_0^\infty t^{p-1}\,|\{x\in\Omega ;\ |u(x)|>t\}|\,dt$$
has  no analogue in the framework of variable exponents.

(ii) Spaces $L^{p(x)}$ do {\it not} satisfy the {\it mean continuity property}. More exactly, if $p$ is nonconstant and continuous in an open ball $B$, then there is some $u\in L^{p(x)}(B)$ such that $u(x+h)\not\in L^{p(x)}(B)$ for every $h\in\RR^N$ with arbitrary small norm.

(iii) Function spaces with variable exponent
 are {\it never} invariant with respect to translations.  The convolution is also limited. For instance,  the classical Young inequality
$$| f*g|_{p(x)}\leq C\, | f|_{p(x)}\, \| g\|_{L^1}$$
remains true if and only if
$p$ is constant.

\textbf{Proposition 2.6} (see \cite{e13}). If the assumption (p$_{1}$) is
satisfied, then $\lambda _{p(\cdot )}$ defined in (\ref{d2}) is positive.
\newline

Next, we prove some results related to the $p(x)$-Laplace operator $-\Delta
_{p(x)}$ as defined at the beginning of Section 1. Consider the following
functional
\begin{equation*}
J(u)=\int_{\Omega }\frac{1}{p(x)}\left\vert \nabla u\right\vert ^{p(x)}dx,%
\text{ }u\in X:=W_{0}^{1,p(\cdot )}(\Omega ).
\end{equation*}

Then (see \cite{e39}) $J\in C^{1}(X, \mathbb{R})$ and the $p(x)$-Laplace
operator is the derivative operator of $J$ in the weak sense. We denote $%
L=J^{^{\prime }}:X\rightarrow X^{\ast }$, then
\begin{equation*}
(L(u),v)=\int_{\Omega }\left\vert \nabla u\right\vert ^{p(x)-2}\nabla
u\nabla vdx\text{, }\forall v,u\in X.
\end{equation*}

\textbf{Theorem 2.7 }(see \cite{e12,e17}). i) $L:X\rightarrow X^{\ast }$ is
a continuous, bounded and strictly monotone operator;

ii) $L$ is a mapping of type $(S_{+})$, that is, if $u_{n}\rightharpoonup u$
in $X$ and $\underset{n\rightarrow +\infty }{\overline{\lim }}$ $%
(L(u_{n})-L(u),u_{n}-u)\leq 0,$ then $u_{n}\rightarrow u$ in $X$;

iii) $L:X\rightarrow X^{*}$ is a homeomorphism.\newline

\medskip Denote
\begin{equation*}
B(x_{0},\varepsilon ,\delta ,\theta )=\left\{x\in \mathbb{R} ^{N}\mid \delta
\leq \left\vert x-x_{0}\right\vert \leq \varepsilon ,\frac{x-x_{0}}{
\left\vert x-x_{0}\right\vert }\cdot \frac{\nabla p(x_{0})}{\left\vert
\nabla p(x_{0})\right\vert }\geq \cos \theta \right\},
\end{equation*}
where $\theta \in (0,\frac{ \pi }{2}).$ Then we obtain the following.\newline

\textbf{Lemma 2.8}. If $p\in C^{1}(\overline{\Omega })$, $x_{0}\in \Omega $
satisfy $\nabla p(x_{0})\neq 0$, then there exists small enough $\varepsilon
>0$ such that
\begin{equation}
(x-x_{0})\cdot \nabla p(x)>0,\ \forall x\in B(x_{0},\varepsilon ,\delta
,\theta )\text{,}  \label{a3}
\end{equation}%
and
\begin{equation}
\max \{p(x)\mid x\in \overline{B(x_{0},\varepsilon )}\}=\max \{p(x)\mid x\in
B(x_{0},\varepsilon ,\delta ,\theta ),\left\vert x-x_{0}\right\vert
=\varepsilon \}.  \label{a4}
\end{equation}

\textbf{Proof.} A proof of this lemma can be found in \cite{e32}. For
readers' convenience, we include it here.

Since $p\in C^{1}(\overline{\Omega })$, for any $x\in B(x_{0},\varepsilon
,\delta ,\theta )$, when $\varepsilon >0$ is small enough, we have
\begin{eqnarray*}
\nabla p(x)\cdot (x-x_{0}) &=&(\nabla p(x_{0})+o(1))\cdot (x-x_{0}) \\
&=&\nabla p(x_{0})\cdot (x-x_{0})+o(\left\vert x-x_{0}\right\vert ) \\
&\geq &\left\vert \nabla p(x_{0})\right\vert \left\vert x-x_{0}\right\vert
\cos \theta +o(\left\vert x-x_{0}\right\vert )>0,
\end{eqnarray*}%
where $o(1)\in
%TCIMACRO{\U{211d} }%
%BeginExpansion
\mathbb{R}
%EndExpansion
^{N}$ is a function and $o(1)\rightarrow 0$ uniformly as $\left\vert
x-x_{0}\right\vert \rightarrow 0$.

When $\varepsilon $ is small enough, condition (\ref{a3}) is valid. Since $%
p\in C^{1}(\overline{\Omega })$, there exist a small enough positive $%
\varepsilon $ such that%
\begin{equation*}
p(x)-p(x_{0})=\nabla p(y)\cdot (x-x_{0})=(\nabla p(x_{0})+o(1))\cdot
(x-x_{0}),\text{ }
\end{equation*}%
where $y=x_{0}+\tau (x-x_{0})$ and $\tau \in (0,1)$, $o(1)\in
%TCIMACRO{\U{211d} }%
%BeginExpansion
\mathbb{R}
%EndExpansion
^{N}$ is a function and $o(1)\rightarrow 0$ uniformly as $\left\vert
x-x_{0}\right\vert \rightarrow 0$.

Suppose that $x\in \overline{B(x_{0},\varepsilon )}\backslash
B(x_{0},\varepsilon ,\delta ,\theta )$. Denote $x^{\ast }=x_{0}+\varepsilon
\nabla p(x_{0})/\left\vert \nabla p(x_{0})\right\vert $.

Suppose that $\frac{x-x_{0}}{\left\vert x-x_{0}\right\vert }\cdot \frac{%
\nabla p(x_{0})}{\left\vert \nabla p(x_{0})\right\vert }<\cos \theta $. When
$\varepsilon $ is small enough, we have
\begin{eqnarray*}
p(x)-p(x_{0}) &=&(\nabla p(x_{0})+o(1))\cdot (x-x_{0}) \\
&<&\left\vert \nabla p(x_{0})\right\vert \left\vert x-x_{0}\right\vert \cos
\theta +\varepsilon \cdot o(1) \\
&\leq &(\nabla p(x_{0})+o(1))\cdot \varepsilon \nabla p(x_{0})/\left\vert
\nabla p(x_{0})\right\vert  \\
&=&p(x^{\ast })-p(x_{0}),
\end{eqnarray*}%
where $o(1)\in
%TCIMACRO{\U{211d} }%
%BeginExpansion
\mathbb{R}
%EndExpansion
^{N}$ is a function and $o(1)\rightarrow 0$ as $\varepsilon \rightarrow 0$.

Suppose that $\left\vert x-x_{0}\right\vert <\delta $. When $\varepsilon $
is small enough, we have%
\begin{eqnarray*}
p(x)-p(x_{0}) &=&(\nabla p(x_{0})+o(1))\cdot (x-x_{0}) \\
&\leq &\left\vert \nabla p(x_{0})\right\vert \left\vert x-x_{0}\right\vert
+\varepsilon \cdot o(1) \\
&<&(\nabla p(x_{0})+o(1))\cdot \varepsilon \nabla p(x_{0})/\left\vert \nabla
p(x_{0})\right\vert  \\
&=&p(x^{\ast })-p(x_{0}),
\end{eqnarray*}%
where $o(1)\in
%TCIMACRO{\U{211d} }%
%BeginExpansion
\mathbb{R}
%EndExpansion
^{N}$ is a function and $o(1)\rightarrow 0$ as $\varepsilon \rightarrow 0$.
Thus
\begin{equation}
\max \{p(x)\mid x\in \overline{B(x_{0},\varepsilon )}\}=\max \{p(x)\mid x\in
B(x_{0},\varepsilon ,\delta ,\theta )\}.  \label{d1}
\end{equation}

It follows from (\ref{a3}) and (\ref{d1}) that relation (\ref{a4}) holds.

The proof of Lemma 2.8 is thus complete. $\square $\newline

\textbf{Lemma 2.9} Suppose that $F(x,u)$ satisfies (f$_{4}$). Let
\begin{equation*}
h(x)=\left\{
\begin{array}{cc}
0, & \quad\mbox{if}\ \left\vert x-x_{0}\right\vert >\varepsilon \\
\varepsilon -\left\vert x-x_{0}\right\vert , & \quad\mbox{if}\ \left\vert
x-x_{0}\right\vert \leq \varepsilon%
\end{array}%
\right. ,
\end{equation*}%
where $\varepsilon $ is defined as in Lemma 2.8. Then there exists large
enough $t$ such that
\begin{equation*}
\int_{\Omega }\left\vert \nabla th\right\vert ^{p(x)}dx-\int_{\Omega
}F(x,th)dx\rightarrow -\infty \text{ as }t\rightarrow +\infty .
\end{equation*}

\textbf{Proof.} Obviously,
\begin{equation*}
\int_{\Omega }\frac{1}{p(x)}\left\vert \nabla th\right\vert ^{p(x)}dx\leq
C_{2}\int_{B(x_{0},\varepsilon ,\delta ,\theta )}\left\vert \nabla
th\right\vert ^{p(x)}dx.
\end{equation*}

We make a spherical coordinate transformation. Denote $r=\left\vert
x-x_{0}\right\vert $. Since $p\in C^{1}(\overline{\Omega })$, it follows
from (\ref{a3}) that there exist positive constants $c_{1}$ and $c_{2}$ such
that

\begin{equation*}
p(\varepsilon ,\omega )-c_{2}(\varepsilon -r)\leq p(r,\omega )\leq
p(\varepsilon ,\omega )-c_{1}(\varepsilon -r),\forall (r,\omega )\in
B(x_{0},\varepsilon ,\delta ,\theta ).
\end{equation*}

Therefore
\begin{eqnarray}
\int_{B(x_{0},\varepsilon ,\delta ,\theta )}\left\vert \nabla th\right\vert
^{p(x)}dx &=&\int_{B(x_{0},\varepsilon ,\delta ,\theta )}\left\vert
t\right\vert ^{p(r,\omega )}r^{N-1}drd\omega  \notag \\
&\leq &\int_{B(x_{0},\varepsilon ,\delta ,\theta )}\left\vert t\right\vert
^{p(\varepsilon ,\omega )-c_{1}(\varepsilon -r)}r^{N-1}drd\omega  \notag \\
&\leq &\varepsilon ^{N-1}\int_{B(x_{0},\varepsilon ,\delta ,\theta
)}t^{p(\varepsilon ,\omega )-c_{1}(\varepsilon -r)}drd\omega  \notag \\
&\leq &\varepsilon ^{N-1}\int_{B(x_{0},1,1,\theta )}\frac{t^{p(\varepsilon
,\omega )}}{c_{1}\ln t}d\omega .  \label{a5}
\end{eqnarray}

Denote
\begin{equation*}
G(x,u)=\frac{F(x,u)}{\left\vert u\right\vert ^{p(x)}[\ln (e+\left\vert
u\right\vert )]^{p(x)}}.
\end{equation*}

Then
\begin{equation}
G(x,u)\rightarrow +\infty \text{ uniformly as }\left\vert u\right\vert
\rightarrow +\infty \text{ for }x\in \overline{\Omega }.  \label{a8}
\end{equation}

Thus there exists a positive constant $M$ such that
\begin{equation*}
G(x,u)\geq 1,\forall \left\vert u\right\vert \geq M,\forall x\in \overline{%
\Omega }.
\end{equation*}

Denote
\begin{eqnarray*}
E_{1} &=&\left\{ x\in B(x_{0},\varepsilon )\mid th\geq M\right\} =\left\{
x\in B(x_{0},\varepsilon )\mid \left\vert x-x_{0}\right\vert \leq
\varepsilon -\frac{M}{t}\right\} , \\
E_{2} &=&B(x_{0},\varepsilon )\backslash E_{1}.
\end{eqnarray*}

Then we have%
\begin{eqnarray*}
\int_{\Omega }F(x,th)dx &=&\int_{B(x_{0},\varepsilon )}F(x,th)dx \\
&=&\int_{E_{1}}F(x,th)dx+\int_{E_{2}}F(x,th)dx \\
&\geq &\int_{E_{1}}F(x,th)dx-C_{1}.
\end{eqnarray*}

When $t$ is large enough, we have%
\begin{eqnarray*}
&&\int_{E_{1}}F(x,th)dx \\
&=&\int_{E_{1}}\left\vert th\right\vert ^{p(x)}[\ln (e+\left\vert
th\right\vert )]^{p(x)}G(x,th)dx \\
&=&\int_{B(x_{0},\varepsilon -\frac{M}{t},\delta ,\theta )}C_{1}\left\vert
th\right\vert ^{p(x)}[\ln (e+\left\vert th\right\vert )]^{p(x)}G(x,th)dx \\
&=&\int_{B(x_{0},\varepsilon -\frac{M}{t},\delta ,\theta )}C_{1}\left\vert
t(\varepsilon -r)\right\vert ^{p(r,\omega )}r^{N-1}[\ln (e+\left\vert
t(\varepsilon -r)\right\vert )]^{p(r,\omega )}G(r,\omega ,t(\varepsilon
-r))drd\omega \\
&\geq &C_{1}\delta ^{N-1}\int_{B(x_{0},\varepsilon -\frac{M}{t},\delta
,\theta )}\left\vert t\right\vert ^{p(\varepsilon ,\omega
)-c_{2}(\varepsilon -r)}\left\vert \varepsilon -r\right\vert ^{p(\varepsilon
,\omega )-c_{1}(\varepsilon -r)} \\
&&[\ln (e+\left\vert t(\varepsilon -r)\right\vert )]^{p(r,\omega
)}G(r,\omega ,t(\varepsilon -r))drd\omega \\
&=&C_{1}\delta ^{N-1}\int_{B(x_{0},1,1,\theta )}d\omega\\
&& \int_{\delta
}^{\varepsilon -\frac{M}{t}}\left\vert t\right\vert ^{p(\varepsilon ,\omega
)-c_{2}(\varepsilon -r)}\left\vert \varepsilon -r\right\vert ^{p(\varepsilon
,\omega )-c_{1}(\varepsilon -r)} 
[\ln (e+\left\vert t(\varepsilon -r)\right\vert )]^{p(r,\omega
)}G(r,\omega ,t(\varepsilon -r))dr \\
&\geq &C_{1}\delta ^{N-1}\int_{B(x_{0},1,1,\theta )}d\omega  \\
&&
\int_{\delta
}^{\varepsilon -\frac{1}{\ln t}}\left\vert t\right\vert ^{p(\varepsilon
,\omega )-c_{2}(\varepsilon -r)}\left\vert \varepsilon -r\right\vert
^{p(\varepsilon ,\omega )}
[\ln (e+\left\vert t(\varepsilon -r)\right\vert )]^{p(r,\omega
)}G(r,\omega ,t(\varepsilon -r))dr \\
&\geq &C_{2}\delta ^{N-1}G(r_{t},\omega _{t},t(\varepsilon
-r_{t}))\\
&&\int_{B(x_{0},1,1,\theta )}(\frac{1}{\ln t})^{p(\varepsilon ,\omega
)}[\ln (e+\frac{t}{\ln t})]^{p(\varepsilon ,\omega )}
\int_{\delta
}^{\varepsilon -\frac{1}{\ln t}}\left\vert t\right\vert ^{p(\varepsilon
,\omega )-c_{2}(\varepsilon -r)}drd\omega \\
&\geq &C_{3}\delta ^{N-1}G(r_{t},\omega _{t},t(\varepsilon
-r_{t}))\int_{B(x_{0},1,1,\theta )}\frac{\left\vert t\right\vert
^{p(\varepsilon ,\omega )-\frac{c_{2}}{\ln t}}}{c_{2}\ln t}d\omega \\
&\geq &C_{4}\delta ^{N-1}G(r_{t},\omega _{t},t(\varepsilon
-r_{t}))\int_{B(x_{0},1,1,\theta )}\frac{\left\vert t\right\vert
^{p(\varepsilon ,\omega )}}{c_{2}\ln t}d\omega ,
\end{eqnarray*}%
where $(r_{t},\omega _{t})\in E_{1}$ is such that
\begin{equation*}
G(r_{t},\omega _{t},t(\varepsilon -r_{t}))=\min \left\{ G(r,\omega
,t(\varepsilon -r))\mid (r,\omega )\in B(x_{0},\varepsilon -\frac{1}{\ln t}%
,\delta ,\theta )\right\} .
\end{equation*}

Note that $t(\varepsilon -r_{t})\geq \frac{t}{\ln t}\rightarrow +\infty $ as
$t\rightarrow +\infty $. Thus
\begin{equation}
\int_{\Omega }F(x,th)dx\geq G(r_{t},\omega _{t},t(\varepsilon
-r_{t}))C_{5}\int_{B(x_{0},1,1,\theta )}\frac{\left\vert t\right\vert
^{p(\varepsilon ,\omega )}}{\ln t}d\omega -C_{1}\text{ as }t\rightarrow
+\infty .  \label{a6}
\end{equation}

It follows from (\ref{a5}), (\ref{a8}) and (\ref{a6}) that $\Psi
(th)\rightarrow -\infty $. Proof of Lemma 2.9 is thus complete. $\square $

\textbf{Lemma 2.10} The following $K_{i}$ ($i=1,2,3$) satisfy the hypothesis
($K$)
\begin{eqnarray*}
K_{1}(t) &=&\ln (e+\left\vert t\right\vert ); \\
K_{2}(t) &=&\ln (e+\ln (e+\left\vert t\right\vert ));\ \mbox{and} \\
K_{3}(t) &=&[\ln (e+\ln (e+\left\vert t\right\vert ))]\ln (e+\left\vert
t\right\vert ).
\end{eqnarray*}

\textbf{Proof.} We only need to check that $K_{3}(t)$ satisfies the
hypothesis ($K$). The proofs for the other functions are similar.

We observe that $1\leq K(\cdot )\in C^{1}([0,+\infty ),[1,+\infty ))$ is
increasing and $K(t)\rightarrow +\infty $ as $t\rightarrow +\infty $. So we
only need to prove that $tK^{\prime }(t)/K(t)\leq \sigma \in (0,1)$, where $%
\sigma $ is a constant. By computation we obtain
\begin{eqnarray*}
\frac{tK^{\prime }}{K} &=&\frac{t}{K}\left\{\frac{[\ln (e+\left\vert
t\right\vert )]\mathrm{sgn}\,t}{[e+\ln (e+\left\vert t\right\vert
)](e+\left\vert t\right\vert )}+\frac{[\ln (e+\ln (e+\left\vert t\right\vert
))]\mathrm{sgn}\,t}{(e+\left\vert t\right\vert )}\right\} \\
&=&\frac{\left\vert t\right\vert }{[\ln (e+\ln (e+\left\vert t\right\vert
))][e+\ln (e+\left\vert t\right\vert )](e+\left\vert t\right\vert )} \\
&&+\frac{\left\vert t\right\vert }{[\ln (e+\left\vert t\right\vert
)](e+\left\vert t\right\vert )}\,.
\end{eqnarray*}

We have
\begin{eqnarray*}
\left\vert t\right\vert &\leq &\frac{1}{3}[\ln (e+\ln (e+\left\vert
t\right\vert ))][e+\ln (e+\left\vert t\right\vert )](e+\left\vert
t\right\vert ), \\
\left\vert t\right\vert &\leq &\frac{1}{2}[\ln (e+\left\vert t\right\vert
)](e+\left\vert t\right\vert )
\end{eqnarray*}
and we complete the proof by observing that
\begin{equation*}
\frac{tK^{\prime }}{K}\leq \frac{5}{6},\ \forall t\in
%TCIMACRO{\U{211d} }%
%BeginExpansion
\mathbb{R}
%EndExpansion
.
\end{equation*}

\section{ Proofs of main results}

In this section we give the proofs of our main results.\newline

\textbf{Definition 3.1 }We say that $u\in W_{0}^{1,p(\cdot )}(\Omega )$ is a
weak solution of (P) if
\begin{equation*}
\int_{\Omega }\left\vert \nabla u\right\vert ^{p(x)-2}\nabla u\cdot \nabla
vdx=\int_{\Omega }f(x,u)vdx,\quad\forall v\in X:=W_{0}^{1,p(\cdot )}(\Omega ).
\end{equation*}

The corresponding functional of (P) is
\begin{equation*}
\varphi \left( u\right) =\int_{\Omega }\frac{1}{p(x)}\left\vert \nabla
u\right\vert ^{p(x)}dx-\int_{\Omega }F(x,u)dx,\quad\forall u\in X,
\end{equation*}%
where $F(x,t)=\int_{0}^{t}f(x,s)ds$.\newline

\textbf{Definition 3.2 }We say that $\varphi $ satisfies the Cerami
condition in $X$, if any sequence $\left\{ u_{n}\right\} \subset X$ such
that $\left\{ \varphi (u_{n})\right\} $ is bounded and $\left\Vert \varphi
^{\prime }(u_{n})\right\Vert (1+\left\Vert u_{n}\right\Vert )\rightarrow 0$
as $n\rightarrow +\infty $ has a convergent subsequence.\newline

\textbf{Lemma 3.3} If $f$ satisfies (f$_{0}$) and (f$_{1}$), then $\varphi $
satisfies the Cerami condition.\newline

\textbf{Proof}. Let $\{u_{n}\}\subset X$ be a Cerami sequence, that is $%
\varphi (u_{n})\rightarrow c$ and $\left\Vert \varphi ^{\prime
}(u_{n})\right\Vert (1+\left\Vert u_{n}\right\Vert )\rightarrow 0$.
Therefore $\varphi ^{\prime }(u_{n})=L(u_{n})-f(x,u_{n})\rightarrow 0$ in $%
X^{\ast }$, then we have $L(u_{n})=f(x,u_{n})+o_{n}(1)$, where $%
o_{n}(1)\rightarrow 0$ in $X^{\ast }$ as $n\rightarrow \infty $. Suppose
that $\{u_{n}\}$ is bounded, then $\{u_{n}\}$ has a weakly convergent
subsequence in $X$. Without loss of generality, we assume that $%
u_{n}\rightharpoonup u$, then by Proposition 2.2 and 2.5, we have $%
f(x,u_{n})\rightarrow f(x,u)$ in $X^{\ast }$. Thus $%
L(u_{n})=f(x,u_{n})+o_{n}(1)\rightarrow f(x,u)$ in $X^{\ast }$. Since $L$ is
a homeomorphism, we have $u_{n}\rightarrow L^{-1}(f(x,u))$ in $X$, and so $%
\varphi $ satisfies the Cerami condition. Therefore $u=L^{-1}(f(x,u))$, then
$L(u)=f(x,u)$, this means $u$ is a solution of (P). Thus we only need to
prove the boundedness of the Cerami sequence $\{u_{n}\}$.

We argue by contradiction. Then there exist $c\in \mathbb{R}$ and $%
\{u_{n}\}\subset X$ satisfying:%
\begin{equation*}
\varphi (u_{n})\rightarrow c\text{, }\left\Vert \varphi ^{\prime
}(u_{n})\right\Vert (1+\left\Vert u_{n}\right\Vert )\rightarrow 0\text{, }%
\left\Vert u_{n}\right\Vert \rightarrow +\infty .
\end{equation*}

Obviously,%
\begin{equation*}
\left\vert \frac{1}{p(x)}u_{n}\right\vert _{p(\cdot )}\leq \frac{1}{p^{-}}%
\left\vert u_{n}\right\vert _{p(\cdot )},\left\vert \nabla \frac{1}{p(x)}%
u_{n}\right\vert _{p(\cdot )}\leq \frac{1}{p^{-}}\left\vert \nabla
u_{n}\right\vert _{p(\cdot )}+C\left\vert u_{n}\right\vert _{p(\cdot )}.
\end{equation*}

Thus $\left\Vert \frac{1}{p(x)}u_{n}\right\Vert \leq C\left\Vert
u_{n}\right\Vert $. Therefore $(\varphi ^{\prime }(u_{n}),\frac{1}{p(x)}%
u_{n})\rightarrow 0$. We may assume that%
\begin{eqnarray*}
c+1 &\geq &\varphi (u_{n})-(\varphi ^{\prime }(u_{n}),\frac{1}{p(x)}u_{n}) \\
&=&\int_{\Omega }\frac{1}{p(x)}\left\vert \nabla u_{n}\right\vert
^{p(x)}dx-\int_{\Omega }F(x,u_{n})dx \\
&&-\{\int_{\Omega }\frac{1}{p(x)}\left\vert \nabla u_{n}\right\vert
^{p(x)}dx-\int_{\Omega }\frac{1}{p(x)}f(x,u_{n})u_{n}dx- \\
&&
\int_{\Omega }\frac{1%
}{p^{2}(x)}u_{n}\left\vert \nabla u_{n}\right\vert ^{p(x)-2}\nabla
u_{n}\nabla pdx\} \\
&\geq &\int_{\Omega }\frac{1}{p^{2}(x)}u_{n}\left\vert \nabla
u_{n}\right\vert ^{p(x)-2}\nabla u_{n}\nabla pdx+\int_{\Omega }\{\frac{1}{%
p(x)}f(x,u_{n})u_{n}-F(x,u_{n})\}dx.
\end{eqnarray*}

Hence
\begin{eqnarray}
\int_{\Omega }\{\frac{f(x,u_{n})u_{n}}{p(x)}-F(x,u_{n})\}dx &\leq
&C_{0}(\int_{\Omega }\left\vert u_{n}\right\vert \left\vert \nabla
u_{n}\right\vert ^{p(x)-1}dx+1)  \notag \\
&\leq &\sigma \int_{\Omega }\frac{\left\vert \nabla u_{n}\right\vert ^{p(x)}%
}{K(\left\vert u_{n}\right\vert )}dx+C_{1}  \notag \\
&&+C(\sigma )\int_{\Omega }\left\vert u_{n}\right\vert ^{p(x)}[K(\left\vert
u_{n}\right\vert )]^{p(x)-1}dx,  \label{a1}
\end{eqnarray}%
where $\sigma $ is a small enough positive constant.

Due to hypothesis (K), it is easy to check that $\frac{u_{n}}{K(\left\vert
u_{n}\right\vert )}\in X$, and $\left\Vert \frac{u_{n}}{K(\left\vert
u_{n}\right\vert )}\right\Vert \leq C_{2}\left\Vert u_{n}\right\Vert $. Let $%
\frac{u_{n}}{K(\left\vert u_{n}\right\vert )}$ be a test function. We have%
\begin{eqnarray*}
&&\int_{\Omega }f(x,u_{n})\frac{u_{n}}{K(\left\vert u_{n}\right\vert )}dx \\
&=&\int_{\Omega }\left\vert \nabla u_{n}\right\vert ^{p(x)-2}\nabla
u_{n}\nabla \frac{u_{n}}{K(\left\vert u_{n}\right\vert )}dx+o(1) \\
&=&\int_{\Omega }\frac{\left\vert \nabla u_{n}\right\vert ^{p(x)}}{%
K(\left\vert u_{n}\right\vert )}dx-\int_{\Omega }u_{n}\left\vert \nabla
u_{n}\right\vert ^{p(x)-2}\nabla u_{n}\nabla \frac{1}{K(\left\vert
u_{n}\right\vert )}dx+o(1).
\end{eqnarray*}

By computation, we obtain
\begin{eqnarray*}
&&\left\vert \int_{\Omega }u_{n}\left\vert \nabla u_{n}\right\vert
^{p(x)-2}\nabla u_{n}\nabla \frac{1}{K(\left\vert u_{n}\right\vert )}%
dx\right\vert \\
&\leq &\int_{\Omega }\left\vert u_{n}\right\vert \left\vert \nabla
u_{n}\right\vert ^{p(x)-1}\frac{\left\vert \nabla K(\left\vert
u_{n}\right\vert )\right\vert }{K^{2}(\left\vert u_{n}\right\vert )}dx \\
&\leq &\int_{\Omega }\frac{\left\vert \nabla u_{n}\right\vert ^{p(x)}}{%
K(\left\vert u_{n}\right\vert )}\frac{\left\vert u_{n}\right\vert K^{\prime
}(\left\vert u_{n}\right\vert )}{K(\left\vert u_{n}\right\vert )}dx.
\end{eqnarray*}

Note that $\frac{\left\vert u_{n}\right\vert K^{\prime }(\left\vert
u_{n}\right\vert )}{K(\left\vert u_{n}\right\vert )}\leq \sigma _{0}\in
(0,1) $. Thus
\begin{equation}
C_{3}\int_{\Omega }\frac{\left\vert \nabla u_{n}\right\vert ^{p(x)}}{%
K(\left\vert u_{n}\right\vert )}dx-C_{4}\leq \int_{\Omega }\frac{%
f(x,u_{n})u_{n}}{K(\left\vert u_{n}\right\vert )}dx\leq C_{5}\int_{\Omega }%
\frac{\left\vert \nabla u_{n}\right\vert ^{p(x)}}{K(\left\vert
u_{n}\right\vert )}dx+C_{6}.  \label{a2}
\end{equation}

By (\ref{a1}), (\ref{a2}) and conditions (f$_{0}$) and (f$_{1}$), we have
\begin{eqnarray*}
&&\int_{\Omega }f(x,u_{n})\frac{u_{n}}{K(\left\vert u_{n}\right\vert )}dx \\
&&\overset{\text{(f}_{1}\text{)}}{\leq }C_{7}\int_{\Omega }\{\frac{%
f(x,u_{n})u_{n}}{p(x)}-F(x,u_{n})\}dx+C_{7} \\
&\leq &C_{7}\{\sigma \int_{\Omega }\frac{\left\vert \nabla u_{n}\right\vert
^{p(x)}}{K(\left\vert u_{n}\right\vert )}dx+C_{8}+C(\sigma )\int_{\Omega
}\left\vert u_{n}\right\vert ^{p(x)}[K(\left\vert u_{n}\right\vert
)]^{p(x)-1}dx\} \\
&\leq &C_{7}\sigma \int_{\Omega }\frac{\left\vert \nabla u_{n}\right\vert
^{p(x)}}{K(\left\vert u_{n}\right\vert )}dx+C_{7}C(\sigma )\int_{\Omega
}\left\vert u_{n}\right\vert ^{p(x)}[K(\left\vert u_{n}\right\vert
)]^{p(x)-1}dx+C_{9} \\
&&\overset{(\ref{a2})}{\leq }\frac{1}{2}\int_{\Omega }\frac{f(x,u_{n})u_{n}}{%
K(\left\vert u_{n}\right\vert )}dx+C_{7}C(\sigma )\int_{\Omega }\left\vert
u_{n}\right\vert ^{p(x)}[K(\left\vert u_{n}\right\vert )]^{p(x)-1}dx+C_{10}.
\end{eqnarray*}

Thus, by condition (f$_{1}$) and the above inequality, we can see
\begin{eqnarray}
&&\int_{\Omega }f(x,u_{n})\frac{u_{n}}{K(\left\vert u_{n}\right\vert )}dx
\notag \\
&\leq &C_{11}\int_{\Omega }\left\vert u_{n}\right\vert ^{p(x)}[K(\left\vert
u_{n}\right\vert )]^{p(x)-1}dx+C_{12}.  \label{b1}
\end{eqnarray}%
Note that $\frac{tf(x,t)}{\left\vert t\right\vert ^{p(x)}[K(t)]^{p(x)}}%
\rightarrow +\infty $ uniformly as $\left\vert t\right\vert \rightarrow
+\infty $ for $x\in \overline{\Omega }$. We claim that
\begin{equation*}
\int_{\Omega }\left\vert u_{n}\right\vert ^{p(x)}[K(\left\vert
u_{n}\right\vert )]^{p(x)-1}dx\quad \mbox{is
bounded.}
\end{equation*}%
This means that
\begin{equation*}
\int_{\Omega }f(x,u_{n})\frac{u_{n}}{K(\left\vert u_{n}\right\vert )}dx\quad %
\mbox{is bounded.}
\end{equation*}

In fact, by (K), we observe that there exists $M>0$ large enough such that
\begin{equation}
\frac{tf(x,t)}{K(t)}>2C_{11}\left\vert t\right\vert
^{p(x)}[K(t)]^{p(x)-1},\forall \left\vert t\right\vert \geq M.  \label{b2}
\end{equation}

Denote $\Omega _{n}=\{x\in \Omega \mid \left\vert u_{n}\right\vert \geq M\}$%
. We have
\begin{equation}
\int_{\Omega }f(x,u_{n})\frac{u_{n}}{K(\left\vert u_{n}\right\vert )}dx\geq
\int_{\Omega _{n}}2C_{11}\left\vert u_{n}\right\vert ^{p(x)}[K(\left\vert
u_{n}\right\vert )]^{p(x)-1}dx-C_{12}.  \label{b3}
\end{equation}

Combining (\ref{b1})-(\ref{b3}), we obtain
\begin{equation*}
\int_{\Omega _{n}}C_{11}\left\vert u_{n}\right\vert ^{p(x)}[K(\left\vert
u_{n}\right\vert )]^{p(x)-1}dx\leq C_{13},
\end{equation*}%
and hence
\begin{equation*}
\int_{\Omega }C_{11}\left\vert u_{n}\right\vert ^{p(x)}[K(\left\vert
u_{n}\right\vert )]^{p(x)-1}dx\leq C_{14}.
\end{equation*}

Thus
\begin{equation*}
\int_{\Omega }f(x,u_{n})\frac{u_{n}}{K(\left\vert u_{n}\right\vert )}dx\leq
C_{14},\text{ for any }n=1,2,\cdots .
\end{equation*}

This combine (f$_{0}$) implies that
\begin{equation}
\{\int_{\Omega }\frac{\left\vert f(x,u_{n})u_{n}\right\vert }{K(\left\vert
u_{n}\right\vert )}dx\}\text{ is bounded.}  \label{w1}
\end{equation}

Let $\varepsilon >0$ satisfy $\varepsilon <\min \{1,p^{-}-1,\frac{1}{p^{\ast
+}},(\frac{p^{\ast }}{\alpha })^{-}-1\}$. Since $\left\Vert \varphi ^{\prime
}(u_{n})\right\Vert \left\Vert u_{n}\right\Vert \rightarrow 0$, we get
\begin{eqnarray*}
\int_{\Omega }\left\vert \nabla u_{n}\right\vert ^{p(x)}dx &=&\int_{\Omega
}f(x,u_{n})u_{n}dx+o(1) \\
&\leq &\int_{\Omega }^{\varepsilon }\left\vert f(x,u_{n})u_{n}\right\vert
^{\varepsilon }[K(\left\vert u_{n}\right\vert )]^{1-\varepsilon }\left[
\frac{\left\vert f(x,u_{n})u_{n}\right\vert }{K(\left\vert u_{n}\right\vert )%
}\right] ^{1-\varepsilon }dx+o(1).
\end{eqnarray*}

By condition (f$_{1}$), we have%
\begin{equation*}
\left\vert f(x,u_{n})u_{n}\right\vert \geq \left\vert u_{n}\right\vert
^{p(x)}\text{ for large enough }\left\vert u_{n}\right\vert \text{,}
\end{equation*}%
and
\begin{equation*}
\lbrack K(\left\vert u_{n}\right\vert )]^{1-\varepsilon }\leq \lbrack \ln
(e+\left\vert u_{n}\right\vert )]^{2(1-\varepsilon )}\text{ for large enough
}\left\vert u_{n}\right\vert \text{,}
\end{equation*}%
then we have
\begin{equation*}
\left\vert f(x,u_{n})u_{n}\right\vert ^{\varepsilon }[K(\left\vert
u_{n}\right\vert )]^{1-\varepsilon }\leq C_{15}(\left\vert
f(x,u_{n})u_{n}\right\vert ^{\varepsilon (1+\varepsilon )}+1).
\end{equation*}

Therefore
\begin{eqnarray*}
\int_{\Omega }\left\vert \nabla u_{n}\right\vert ^{p(x)}dx &=&\int_{\Omega
}f(x,u_{n})u_{n}dx+o(1) \\
&\leq &C_{15}(1+\left\Vert u_{n}\right\Vert )^{1+\varepsilon }\int_{\Omega }
\left[ \frac{\left\vert f(x,u_{n})u_{n}\right\vert ^{1+\varepsilon }+1}{%
(1+\left\Vert u_{n}\right\Vert )^{\frac{1+\varepsilon }{\varepsilon }}}%
\right] ^{\varepsilon }\left[ \frac{\left\vert f(x,u_{n})u_{n}\right\vert }{%
K(\left\vert u_{n}\right\vert )}\right] ^{1-\varepsilon }dx+o(1).
\end{eqnarray*}

By Young's inequality, we have
\begin{equation}
\int_{\Omega }\left\vert \nabla u_{n}\right\vert ^{p(x)}dx\leq
C_{15}(1+\left\Vert u_{n}\right\Vert )^{1+\varepsilon }\int_{\Omega }\frac{%
\left\vert f(x,u_{n})u_{n}\right\vert ^{1+\varepsilon }+1}{(1+\left\Vert
u_{n}\right\Vert )^{\frac{1+\varepsilon }{\varepsilon }}}+\frac{\left\vert
f(x,u_{n})u_{n}\right\vert }{K(\left\vert u_{n}\right\vert )}dx+o(1).
\label{w2}
\end{equation}

According to the definition of $\varepsilon $, we have
\begin{equation*}
\left\vert f(x,u_{n})u_{n}\right\vert ^{1+\varepsilon }+1\leq C(\left\vert
u_{n}\right\vert ^{p^{\ast }(x)}+1)
\end{equation*}%
and
\begin{equation*}
(1+\left\Vert u_{n}\right\Vert )^{\frac{1+\varepsilon }{\varepsilon }}\geq
(1+\left\Vert u_{n}\right\Vert )^{(1+\varepsilon )(p^{\ast })^{+}}.
\end{equation*}

Therefore%
\begin{equation*}
\int_{\Omega }\frac{\left\vert f(x,u_{n})u_{n}\right\vert ^{1+\varepsilon }+1%
}{(1+\left\Vert u_{n}\right\Vert )^{\frac{1+\varepsilon }{\varepsilon }}}%
dx\leq \int_{\Omega }\frac{C(\left\vert u_{n}\right\vert ^{p^{\ast }(x)}+1)}{%
(1+\left\Vert u_{n}\right\Vert )^{\frac{1+\varepsilon }{\varepsilon }}}%
dx\leq
\end{equation*}
\begin{equation*}
 \leq\frac{C(\left\vert u_{n}\right\vert ^{(p^{\ast })^{+}}+1)}{%
(1+\left\Vert u_{n}\right\Vert )^{\frac{1+\varepsilon }{\varepsilon }}}\leq
\frac{C_{\#}(\left\Vert u_{n}\right\Vert ^{(p^{\ast })^{+}}+1)}{%
(1+\left\Vert u_{n}\right\Vert )^{\frac{1+\varepsilon }{\varepsilon }}}.
\end{equation*}

Thus, the sequence $$\left\{\di\int_{\Omega }\frac{\left\vert f(x,u_{n})u_{n}\right\vert
^{1+\varepsilon }+1}{(1+\left\Vert u_{n}\right\Vert )^{\frac{1+\varepsilon }{%
\varepsilon }}}dx\right\}$$ is bounded. This combine (\ref{w1}) and (\ref{w2})
implies%
\begin{equation*}
\int_{\Omega }\left\vert \nabla u_{n}\right\vert ^{p(x)}dx\leq
C_{16}(1+\left\Vert u_{n}\right\Vert )^{1+\varepsilon }+C_{17}.
\end{equation*}

Note that $\varepsilon <p^{-}-1$. This is a contradiction, hence $\{u_{n}\}$
is bounded in $X$.

The proof of Lemma 3.3 is thus complete. $\square $

\medskip \textbf{Proof of Theorem 1.1}. We first establish the existence of
a nontrivial weak solution.

We show that $\varphi $ satisfies conditions of the mountain pass lemma. By
Lemma 3.3, $\varphi $ satisfies the Cerami condition. Since $p(x)<\alpha
(x)<p^{\ast }(x)$, the embedding $X\hookrightarrow L^{\alpha (\cdot
)}(\Omega )$ is compact. Hence there exists $C_{0}>0$ such that
\begin{equation*}
\left\vert u\right\vert _{p(\cdot )}\leq C_{0}\left\Vert u\right\Vert \text{%
, }\forall u\in X.
\end{equation*}

Let $\sigma >0$ be small enough such that $\sigma \leq \frac{1}{4}\lambda
_{p(\cdot )}$. By the assumptions (f$_{0}$) and (f$_{2}$), we obtain
\begin{equation*}
F(x,t)\leq \sigma \frac{1}{p(x)}\left\vert t\right\vert ^{p(x)}+C(\sigma
)\left\vert t\right\vert ^{\alpha (x)}\text{, }\forall (x,t)\in \Omega
\times
%TCIMACRO{\U{211d} }%
%BeginExpansion
\mathbb{R}
%EndExpansion
.
\end{equation*}

By (p$_{1}$) and Lemma 2.6, we have $\lambda _{p(\cdot )}>0$ and
\begin{equation*}
\int_{\Omega }\frac{1}{p(x)}\left\vert \nabla u\right\vert ^{p(x)}dx-\sigma
\int_{\Omega }\frac{1}{p(x)}\left\vert u\right\vert ^{p(x)}dx\geq \frac{3}{4}%
\int_{\Omega }\frac{1}{p(x)}\left\vert \nabla u\right\vert ^{p(x)}.
\end{equation*}

Since $\alpha \in C(\overline{\Omega })$ and $p(x)<\alpha (x)<p^{\ast }(x)$,
we can divide the domain $\Omega $ into $n_{0}$ disjoint small subdomains $%
\Omega _{i}$ ($i=1,\cdots ,n_{0}$) such that $\overline{\Omega }=\underset{%
i=1}{\overset{n_{0}}{\cup }}\overline{\Omega _{i}}$ and
\begin{equation*}
\underset{\Omega _{i}}{\sup }\,p(x)<\underset{\Omega _{i}}{\inf }\,\alpha
(x)\leq \underset{\Omega _{i}}{\sup }\,\alpha (x)<\underset{\Omega _{i}}{\inf }\,
p^{\ast }(x).
\end{equation*}

Let
\begin{equation*}
\epsilon =\underset{1\leq i\leq n_{0}}{\min }\{\underset{\Omega _{i}}{\inf }%
\alpha (x)-\underset{\Omega _{i}}{\sup }p(x)\}.
\end{equation*}
and denote by $\left\Vert u\right\Vert _{\Omega _{i}}$ the norm of $u$ on $%
\Omega _{i}$, that is
\begin{equation*}
\int_{\Omega _{i}}\frac{1}{p(x)}\left\vert \nabla \frac{u}{\left\Vert
u\right\Vert _{\Omega _{i}}}\right\vert ^{p(x)}dx+\int_{\Omega _{i}}\frac{1}{%
p(x)}\left\vert \frac{u}{\left\Vert u\right\Vert _{\Omega _{i}}}\right\vert
^{p(x)}dx=1.
\end{equation*}

Then $\left\Vert u\right\Vert _{\Omega _{i}}\leq C\left\Vert u\right\Vert $
and there exist $\xi _{i},\eta _{i}\in \overline{\Omega _{i}}$ such that
\begin{eqnarray*}
\left\vert u\right\vert _{\alpha (\cdot )}^{\alpha (\xi _{i})}
&=&\int_{\Omega _{i}}\left\vert u\right\vert ^{\alpha (x)}dx, \\
\left\Vert u\right\Vert _{\Omega _{i}}^{p(\eta _{i})} &=&\int_{\Omega
_{i}}\left( \frac{1}{p(x)}\left\vert \nabla u\right\vert ^{p(x)}+\frac{1}{%
p(x)}\left\vert u\right\vert ^{p(x)}\right)dx.
\end{eqnarray*}

When $\left\Vert u\right\Vert $ is small enough, we have
\begin{eqnarray*}
C(\sigma )\int_{\Omega }\left\vert u\right\vert ^{\alpha (x)}dx &=&C(\sigma )%
\underset{i=1}{\overset{n_{0}}{\sum }}\int_{\Omega _{i}}\left\vert
u\right\vert ^{\alpha (x)}dx \\
&=&C(\sigma )\underset{i=1}{\overset{n_{0}}{\sum }}\left\vert u\right\vert
_{\alpha (\cdot )}^{\alpha (\xi _{i})}\text{ (where }\xi _{i}\in \overline{%
\Omega _{i}}\text{)} \\
&\leq &C\underset{i=1}{\overset{n_{0}}{\sum }}\left\Vert u\right\Vert
_{\Omega _{i}}^{\alpha (\xi _{i})}\text{ (by Proposition 2.5)} \\
&\leq &C\left\Vert u\right\Vert ^{\epsilon }\underset{i=1}{\overset{n_{0}}{%
\sum }}\left\Vert u\right\Vert _{\Omega _{i}}^{p(\eta _{i})}\text{ (where }%
\eta _{i}\in \overline{\Omega _{i}}\text{)} \\
&=&C\left\Vert u\right\Vert ^{\epsilon }\underset{i=1}{\overset{n_{0}}{\sum }%
}\int_{\Omega _{i}}(\frac{1}{p(x)}\left\vert \nabla u\right\vert ^{p(x)}+%
\frac{1}{p(x)}\left\vert u\right\vert ^{p(x)})dx \\
&=&C\left\Vert u\right\Vert ^{\epsilon }\int_{\Omega }(\frac{1}{p(x)}%
\left\vert \nabla u\right\vert ^{p(x)}+\frac{1}{p(x)}\left\vert u\right\vert
^{p(x)})dx \\
&\leq &\frac{1}{4}\int_{\Omega }\frac{1}{p(x)}\left\vert \nabla u\right\vert
^{p(x)}dx.
\end{eqnarray*}

Thus%
\begin{eqnarray*}
\varphi (u) &\geq &\int_{\Omega }\frac{1}{p(x)}\left\vert \nabla
u\right\vert ^{p(x)}-\sigma \int_{\Omega }\frac{1}{p(x)}\left\vert
u\right\vert ^{p(x)}dx-C(\sigma )\int_{\Omega }\left\vert u\right\vert
^{\alpha (x)}dx \\
&\geq &\frac{1}{2}\int_{\Omega }\frac{1}{p(x)}\left\vert \nabla u\right\vert
^{p(x)}\text{ when }\parallel u\parallel \text{ is small enough.}
\end{eqnarray*}

Therefore, there exist $r>0$ and $\delta >0$ such that $\varphi (u)\geq
\delta >0$ for every $u\in X$ and $\left\Vert u\right\Vert =r$.

Suppose (p$_{2}$) is satisfied. Define $h\in C_{0}\overline{(B(x_{0},3\delta
))}$ as follows:
\begin{equation*}
h(x)=\left\{
\begin{array}{cc}
0, & \left\vert x-x_{0}\right\vert \geq 3\delta \\
3\delta -\left\vert x-x_{0}\right\vert , & 2\delta \leq \left\vert
x-x_{0}\right\vert <3\delta \\
\delta , & \left\vert x-x_{0}\right\vert <2\delta%
\end{array}%
\right. .
\end{equation*}

Note that $$\underset{\left\vert x-x_{0}\right\vert \leq \delta }{\min }p(x)>%
\underset{2\delta \leq \left\vert x-x_{0}\right\vert \leq 3\delta }{\max }%
p(x).$$ It is now easy to check that
\begin{eqnarray*}
\varphi (th) &=&\int_{\Omega }\frac{1}{p(x)}\left\vert \nabla th\right\vert
^{p(x)}-\int_{\Omega }F(x,th)dx \leq \int_{\overline{B(x_{0},3\delta )}\backslash (\overline{%
B(x_{0},2\delta )})}\frac{1}{p(x)}\left\vert \nabla th\right\vert
^{p(x)}-\\
&&\int_{(\overline{B(x_{0},\delta )})}C_{1}\left\vert th\right\vert
^{p(x)}dx+C_{2}\rightarrow -\infty \text{ as }t\rightarrow +\infty .
\end{eqnarray*}

Since $\varphi \left( 0\right) =0,$ $\varphi $ satisfies the conditions of
the mountain pass lemma. So $\varphi $ admits at least one nontrivial
critical point, which implies the problem (P) has a nontrivial weak solution
$u$.

Suppose (f$_{4}$) is satisfied. We may assume that there exists $x_{0}\in
\Omega $ such that $\nabla p(x_{0})\neq 0$.

Define $h\in C_{0}(\overline{B(x_{0},\varepsilon )})$ as follows:
\begin{equation*}
h(x)=\left\{
\begin{array}{cc}
0, & \left\vert x-x_{0}\right\vert \geq \varepsilon \\
\varepsilon -\left\vert x-x_{0}\right\vert , & \left\vert x-x_{0}\right\vert
<\varepsilon%
\end{array}%
\right. .
\end{equation*}

By (f$_{4}$) and Lemma 2.9, there exists $\varepsilon >0$ small enough such
that
\begin{equation*}
\varphi (th)=\int_{\Omega }\frac{1}{p(x)}\left\vert \nabla th\right\vert
^{p(x)}-\int_{\Omega }F(x,th)dx\rightarrow -\infty \text{ as }t\rightarrow
+\infty .
\end{equation*}

Since $\varphi \left( 0\right) =0$, $\varphi $ satisfies the conditions of
the mountain pass lemma. So $\varphi $ admits at least one nontrivial
critical point, which implies that problem (P) has a nontrivial weak
solution $u$. The proof of Theorem 1.1 is thus complete. $\square $

\medskip In order to prove Theorem 1.2, we need to do some preparations.
Note that $X:=W_{0}^{1,p(\cdot )}(\Omega )$ is a reflexive and separable
Banach space (see \cite{22}, Section 17, Theorem 2-3). Therefore there exist
$\left\{ e_{j}\right\} \subset X$ and $\left\{ e_{j}^{\ast }\right\} \subset
X^{\ast } $ such that%
\begin{equation*}
X=\overline{\mathrm{span}}\,\{e_{j}\text{, }j=1,2,\cdots \}\text{, }\left.
{}\right. X^{\ast }=\overline{\mathrm{span}}^{W^{\ast }}\{e_{j}^{\ast }\text{%
, }j=1,2,\cdots \},
\end{equation*}%
and
\begin{equation*}
<e_{j}^{\ast },e_{j}>=\left\{
\begin{array}{c}
1,\quad\mbox{if}\ i=j, \\
0,\quad\mbox{if}\ i\neq j.%
\end{array}%
\right.
\end{equation*}

For convenience, we write $X_{j}=\mathrm{span}\,\{e_{j}\}$, $Y_{k}=\overset{k%
}{\underset{j=1}{\oplus }}X_{j}$ and $Z_{k}=\overline{\overset{\infty }{%
\underset{j=k}{\oplus }}X_{j}}$.\newline

\textbf{Lemma 3.4}. Assume that $\alpha \in C_{+}\left( \overline{\Omega }%
\right) $, $\alpha (x)<p^{\ast }(x)$ for any $x\in \overline{\Omega }$. If
\begin{equation*}
\beta _{k}=\sup \left\{ \left\vert u\right\vert _{\alpha (\cdot )}\left\vert
\left\Vert u\right\Vert =1,u\in Z_{k}\right. \right\} ,
\end{equation*}%
then $\underset{k\rightarrow \infty }{\lim }\beta _{k}=0$.\newline

\textbf{Proof}. Obviously, $0<\beta _{k+1}\leq \beta _{k},$ so $\beta
_{k}\rightarrow \beta \geq 0$. Let $u_{k}\in Z_{k}$ satisfy%
\begin{equation*}
\left\Vert u_{k}\right\Vert =1,0\leq \beta _{k}-\left\vert u_{k}\right\vert
_{\alpha (\cdot )}<\frac{1}{k}.
\end{equation*}%
Then there exists a subsequence of $\{u_{k}\}$ (which we still denote by $%
u_{k}$) such that $u_{k}\rightharpoonup u$, and%
\begin{equation*}
<e_{j}^{\ast },u>=\underset{k\rightarrow \infty }{\lim }\left\langle
e_{j}^{\ast },u_{k}\right\rangle =0\text{, }\forall e_{j}^{\ast }.
\end{equation*}%
This implies that $u=0$, and so $u_{k}\rightharpoonup 0$. Since the
embedding from $W_{0}^{1,p(\cdot )}\left( \Omega \right) $ into $L^{\alpha
(\cdot )}\left( \Omega \right) $ is compact, we can conclude that $%
u_{k}\rightarrow 0$ in $L^{\alpha (\cdot )}\left( \Omega \right) $. Hence we
get $\beta _{k}\rightarrow 0$ as $k\rightarrow \infty $. The proof of Lemma
3.4 is thus complete. $\square $

\medskip In order to prove Theorem 1.2, we need the following auxiliary
result, see \cite[Theorem 4.7]{20a}. If the Cerami condition is replaced by
PS condition, we can use the following property, see \cite[Theorem 3.6]{e39}.

\medskip \textbf{Lemma 3.5}. Suppose that $\varphi \in C^{1}(X, \mathbb{R})$
is even and satisfies the Cerami condition. Let $V^{+}$, $V^{-}\subset X$ be
closed subspaces of $X$ with codim\,$V^{+}+1=$dim $V^{-}$. Suppose that:

($1^{0}$) $\varphi (0)=0$;

($2^{0}$) $\exists \tau >0,$ $\gamma >0$ such that $\forall u\in V^{+}:$ $%
\Vert u\Vert =\gamma \Rightarrow \varphi (u)\geq \tau $; and

($3^{0}$) $\exists \rho >0\ $such that $\forall u\in V^{-}:$ $\Vert u\Vert
\geq \rho \Rightarrow \varphi (u)\leq 0.$

Consider the following set:
\begin{equation*}
\Gamma =\{g\in C^{0}(X,X)\mid g\text{ is odd, }g(u)=u\text{ if }u\in V^{-}%
\text{ and }\Vert u\Vert \geq \rho \}\,.
\end{equation*}%
Then

($a$) $\forall \delta >0$, $g\in \Gamma $, $S_{\delta }^{+}\cap g(V^{-})\neq
\emptyset $, here $S_{\delta }^{+}=\{u\in V^{+}\mid \Vert u\Vert =\delta
\};$ and

($b$) the number $\varpi :=\underset{g\in \Gamma }{\inf }\underset{\text{ }%
u\in V^{-}}{\sup }\varphi (g(u))\geq \tau >0$ is a critical value for $%
\varphi $.\newline

\textbf{Proof of Theorem 1.2}. We first establish the existence of
infinitely many pairs of weak solutions.

According to (f$_{0}$), (f$_{1}$) and (f$_{3}$), $\varphi $ is an even
functional and satisfies the Cerami condition. Let $V_{k}^{+}=Z_{k}$ be a
closed linear subspace of $X$ and $V_{k}^{+}\oplus Y_{k-1}=X$.

Suppose that (f$_{4}$) is satisfied. We may assume that there exists $%
x_{n}\in \Omega $ such that $\nabla p(x_{n})\neq 0$.

Define $h_{n}\in C_{0}(\overline{B(x_{n},\varepsilon _{n})})$ by
\begin{equation*}
h_{n}(x)=\left\{
\begin{array}{cc}
0, & \left\vert x-x_{n}\right\vert \geq \varepsilon _{n} \\
\varepsilon _{n}-\left\vert x-x_{n}\right\vert , & \left\vert
x-x_{n}\right\vert <\varepsilon _{n}%
\end{array}%
\right. .
\end{equation*}

Without loss of generality, we may assume that
\begin{equation*}
\mathrm{supp}\,h_{i}\cap \mathrm{supp}\,h_{j}=\emptyset \text{, }\forall
i\neq j.
\end{equation*}

By Lemma 2.9, we can let $\varepsilon _{n}>0$ be small enough so that
\begin{equation*}
\varphi (th_{n})=\int_{\Omega }\frac{1}{p(x)}\left\vert \nabla
th_{n}\right\vert ^{p(x)}-\int_{\Omega }F(x,th_{n})dx\rightarrow -\infty
\text{ as }t\rightarrow +\infty .
\end{equation*}

Suppose that (p$_{3}$) is satisfied. Define $h_{n}\in C_{0}(\overline{%
B(x_{n},\varepsilon _{n})})$ by
\begin{equation*}
h_{n}(x)=\left\{
\begin{array}{cc}
0, & \left\vert x-x_{n}\right\vert \geq 3\delta _{n} \\
3\delta _{n}-\left\vert x-x_{n}\right\vert , & 2\delta _{n}\leq \left\vert
x-x_{n}\right\vert <3\delta _{n} \\
\delta _{n}, & \left\vert x-x_{n}\right\vert <2\delta _{n}%
\end{array}%
\right. .
\end{equation*}

Note that $\underset{\left\vert x-x_{n}\right\vert \leq \delta _{n}}{\min }%
p(x)>\underset{2\delta _{n}\leq \left\vert x-x_{n}\right\vert \leq 3\delta
_{n}}{\max }p(x)$. It follows that
\begin{eqnarray*}
\varphi (th_{n}) &=&\int_{\Omega }\frac{1}{p(x)}\left\vert \nabla
th_{n}\right\vert ^{p(x)}-\int_{\Omega }F(x,th_{n})dx \\
&\leq &\int_{2\delta _{n}\leq \left\vert x-x_{n}\right\vert \leq 3\delta
_{n}}\frac{1}{p(x)}\left\vert \nabla th_{n}\right\vert
^{p(x)}-\int_{\left\vert x-x_{n}\right\vert \leq \delta _{n}}C_{1}\left\vert
th_{n}\right\vert ^{p(x)}dx+C_{2}\rightarrow -\infty
\end{eqnarray*}
as $t\rightarrow +\infty $.

Set $V_{k}^{-}=\mathrm{span}\,\{h_{1},\cdots ,h_{k}\}$. We will prove that
there exist infinitely many pairs of $V_{k}^{+}$ and $V_{k}^{-}$, such that $%
\varphi $ satisfies the conditions of Lemma 3.5 and the corresponding
critical value satisfies
\begin{equation*}
\varpi _{k}:=\underset{g\in \Gamma }{\inf }\underset{\text{ } u\in V_{k}^{-}}%
{\sup }\varphi (g(u))\rightarrow +\infty
\end{equation*}
when $k\rightarrow +\infty $. This shows that there are infinitely many
pairs of solutions of the problem (P).

For any $m=1,2,\cdots $, we will prove that there exist $\rho _{m}>\gamma
_{m}>0$ and large enough $k_{m}$ such that%
\begin{eqnarray*}
(A_{1})\left. {}\right. \text{ }b_{k_{m}} &:&=\inf \left\{ \varphi (u)\mid
u\in V_{k_{m}}^{+},\left\Vert u\right\Vert =\gamma _{m}\right\} \rightarrow
+\infty \text{ }(m\rightarrow +\infty ); \mbox{and} \\
(A_{2})\left. {}\right. \text{ }a_{k_{m}} &:&=\max \left\{ \varphi
(u)\right\vert \text{ }u\in V_{k_{m}}^{-},\left\Vert u\right\Vert =\rho
_{m}\}\leq 0.
\end{eqnarray*}

First, we prove ($A_{1}$) as follows. By computation, for any $u\in
Z_{k_{m}} $ with $\left\Vert u\right\Vert =\gamma _{m}=m$, we have%
\begin{eqnarray*}
\varphi (u) &=&\int_{\Omega }\frac{1}{p(x)}\left\vert \nabla u\right\vert
^{p(x)}dx-\int_{\Omega }F(x,u)dx \\
&\geq &\frac{1}{p^{+}}\int_{\Omega }\left\vert \nabla u\right\vert
^{p(x)}dx-C\int_{\Omega }\left\vert u\right\vert ^{\alpha
(x)}dx-C_{1}\int_{\Omega }\left\vert u\right\vert dx \\
&\geq &\frac{1}{p^{+}}\left\Vert u\right\Vert ^{p^{-}}-C\left\vert
u\right\vert _{\alpha (\cdot )}^{\alpha (\xi )}-C_{2}\left\vert u\right\vert
_{\alpha (\cdot )}\text{ (where }\xi \in \Omega \text{)} \\
&\geq &\left\{
\begin{array}{l}
\frac{1}{p^{+}}\left\Vert u\right\Vert ^{p^{-}}-C\beta _{k_{m}}^{\alpha
^{-}}\left\Vert u\right\Vert ^{\alpha ^{-}}-C_{2}\beta _{k_{m}}\left\Vert
u\right\Vert \text{, if }\left\vert u\right\vert _{\alpha (\cdot )}\leq 1,
\\
\frac{1}{p^{+}}\left\Vert u\right\Vert ^{p^{-}}-C\beta _{k_{m}}^{\alpha
^{+}}\left\Vert u\right\Vert ^{\alpha ^{+}}-C_{2}\beta _{k_{m}}\left\Vert
u\right\Vert \text{, if }\left\vert u\right\vert _{\alpha (\cdot )}>1,%
\end{array}%
\right. \\
&\geq &\frac{1}{p^{+}}\left\Vert u\right\Vert ^{p^{-}}-C\beta
_{k_{m}}^{\alpha ^{-}}(\left\Vert u\right\Vert ^{\alpha ^{+}}+1)-C_{2}\beta
_{k_{m}}\left\Vert u\right\Vert .
\end{eqnarray*}

Obviously, there exists a large enough $k_{m}$ such that
\begin{equation*}
\frac{1}{p^{+}}\left\Vert u\right\Vert ^{p^{-}}-C\beta _{k_{m}}^{\alpha
^{-}}(\left\Vert u\right\Vert ^{\alpha ^{+}}+1)-C_{2}\beta
_{k_{m}}\left\Vert u\right\Vert \geq \frac{1}{2p^{+}}\left\Vert u\right\Vert
^{p^{-}}\text{, }\forall u\in Z_{k_{m}}\text{ with }\left\Vert u\right\Vert
=\gamma _{m}=m.
\end{equation*}

Therefore $\varphi (u)\geq \frac{1}{2p^{+}}\left\Vert u\right\Vert ^{p^{-}}$%
, $\forall u\in Z_{k_{m}}$ with $\left\Vert u\right\Vert =\gamma _{m}=m$.
Hence $b_{k_{m}}\rightarrow +\infty$ as $m\rightarrow \infty.$

Now we give a proof of ($A_{2}$). According to the above discussion, it is
easy to see that
%\begin{equation*}
$\Psi (th_{k_{m}})\rightarrow -\infty \text{ as }t\rightarrow +\infty .$
%\end{equation*}%
Therefore
\begin{equation*}
\varphi (th)\rightarrow -\infty \text{ as }t\rightarrow +\infty ,\forall
h\in V_{k_{m}}^{-}=\mathrm{span}\,\{h_{1},\cdots ,h_{k_{m}}\}\text{ with}%
\parallel h\parallel =1.
\end{equation*}

This completes the proof of Theorem 1.2. $\square $

\end{document}